\documentclass[12pt,twoside]{article}
\usepackage{amsfonts}
\usepackage{amssymb,amsmath}
\begin{document}
\begin{center} \textbf{A strong ``abc-conjecture"
for certain partitions a+b of c} \vspace{12pt}
\\Constantin M. Petridi
\\ cpetridi@hotmail.com
\end{center}

\begin{center}
\begin{tabular}{p{10cm}}
\begin{small}
\textbf{Abstract} We prove that for any positive integer $c$ and
any $\varepsilon >0$ there are representations of $c$ as a sum
$a+b$ of two coprime positive integers $a$, $b$, such that the
respective radicals $R(abc)$ satisfy
$$ k_{\varepsilon}\;R(c)^{1-\varepsilon}\;c^{2} \hspace{5pt} <
\hspace{5pt} R(abc),$$where $k_{\varepsilon}$ is an absolute
constant depending only on $\varepsilon$. For the representations
in question this is a stronger result than the abc-conjecture
$$\kappa_{\varepsilon}\;c^{\frac{1}{1+\varepsilon}}\hspace{5pt} <
\hspace{5pt} R(abc).$$An upper bound, depending on the number of
prime factors of $c$, is also established.
\end{small}
\end{tabular}
\end{center}

\vspace{10pt}
 \textbf{Preliminaries.} Let
$c=q_{1}^{\alpha_1}q_{2}^{\alpha_2}\cdots
q_{\omega}^{\alpha_\omega},$ \ $q_i$ different primes,
$\alpha_{i}\geq 1,$ be a positive integer. Denote its radical
 $q_{1}q_{2}\cdots q_{\omega},$ by $R(c),$ and similarly $R(n)$ for any integer $n$. Consider the positive
solutions of the Diophantine equation $x + y = c,$ \ $(x,y) = 1,$
\ $x < y $ . Their number is $\varphi(c)/2$. Denoting them in some
order by $a_{i},\hspace{1pt} b_{i}$, $1 \leq i \leq \varphi(c)/2,$
and listing them with their respective radicals, one has
\par
\hspace{1in}
\begin{tabular}{cc}
$a_{1} + b_{1}=c$ & $R(a_{1}b_{1}c)$\\
$a_{2} + b_{2}=c$ & $R(a_{2}b_{2}c)$\\
$ \cdots        $ & $\cdots $        \\
$a_{\frac{\varphi(c)}{2}} + b_{\frac{\varphi(c)}{2}} = c\hspace{23pt}$&$\hspace{35pt} R(a_{\frac{\varphi(c)}{2}}b_{\frac{\varphi(c)}{2}}c).\hspace{10pt}$\\
\end{tabular}\hspace{0.7in}(1)
\par
Form the product of above radicals
\begin{displaymath}
G_{c}=\prod_{1\leq i \leq \frac{\varphi(c)}{2}} R(a_{i}b_{i}c).
\end{displaymath}
$G_{c}^{2/\varphi(c)}$ is the geometric mean of the radicals.
\par
The function $E_{c}(x)$ defined for any real $x\neq 0$ by
\begin{displaymath}
E_{c}(x)=\left[\frac{1}{x}\frac{c}{1}\right]-\sum_{1\leq i \leq
\omega}\left[\frac{1}{x}\frac{c}{q_{i}}\right]+\sum_{1\leq i,j
\leq \omega}\left[\frac{1}{x}\frac{c}{q_{i}q_{j}}\right]-\cdots
+(-1)^{\omega}\left[\frac{1}{x}\frac{c}{q_{1}\cdots
q_{\omega}}\right],
\end{displaymath}
plays a key role in our investigations. If the numbers in the
integral part brackets are integers the function reduces to
\begin{displaymath}
\frac{1}{x}\frac{c}{1}-\sum_{1\leq i \leq
\omega}\frac{1}{x}\frac{c}{q_{i}}+\sum_{1\leq i,j \leq
\omega}\frac{1}{x}\frac{c}{q_{i}q_{j}}-\cdots
+(-1)^{\omega}\frac{1}{x}\frac{c}{q_{1}\cdots
q_{\omega}}=\frac{1}{x}\varphi(c).
\end{displaymath}
Throughout the paper, $p$ designates prime numbers.
\par
\vspace{20pt} \textbf{Theorem 1.}
\begin{displaymath}
G_{c}=R(c)^{\frac{\varphi(c)}{2}}\prod_{\substack{2\leq p < c \\
(p,\;c)=1}}p^{E_{c}(p)}.
\end{displaymath}
\par
\textbf{Proof.} Consider \textit{all} positive solutions of the
Diophantine equation $x + y = c,\quad x \leq y,\hspace{5pt}$ and
their corresponding radicals $R(xyc)$
\par
\begin{tabular}{c}
$1 + (c-1) = c \hspace{1in} R(1(c-1)c)$\\
$2 + (c-2) = c \hspace{1in} R(2(c-2)c)$\\
$\hspace{1.4in} \cdots \hspace{1.9in}\cdots  \hspace{1.2in}$    \\
$\left[\frac{c}{2}\right]+\left(c-\left[\frac{c}{2}\right]\right)=c\hspace{1in}
R\left(\left[\frac{c}{2}\right]\left[c-\frac{c}{2}\right]c\right) $\\
\end{tabular}(2)
\par
Above equalities comprise (1), but include also those for which
$(x,\; y) > 1$.
 That all $q_{i}$ appear $\frac{\varphi(c)}{2}$ times in $G_{c}$ is obvious. Consequently, so does their product
$R(c)$. For primes $p \neq q_{i}$, we apply the
inclusion-exclusion principle. As all numbers $< c$ do occur in
the equalities (2), the number of times $p$ appears in the
radicals (2) is $\left[\frac{c}{p}\right]$. The number of times
$pq_{i}$ appears in the radicals (2) is
$\sum_{i}^{\omega}\left[\frac{c}{pq_{i}} \right]$. The number of
times $pq_{i}q_{j}$ appears in the radicals (2) is $\sum_{i,\; j
}^{\omega}\left[\frac{c}{pq_{i}q_{j}} \right],$ e.t.c. Inserting
these values in the inclusion-exclusion formula we get for the
total number of times p appears in the radicals (1)
\begin{displaymath}
\left[\frac{1}{p}\frac{c}{1}\right]-\sum_{1\leq i \leq
\omega}\left[\frac{1}{p}\frac{c}{q_{i}}\right]+\sum_{1\leq i,j
\leq \omega}\left[\frac{1}{p}\frac{c}{q_{i}q_{j}}\right]-\cdots
+(-1)^{\omega}\left[\frac{1}{p}\frac{c}{q_{1}\cdots
q_{\omega}}\right]=E_{p}(c),
\end{displaymath}
as stated in the theorem.
\par
\vspace{20pt} \textbf{Corollary.} If d denotes the divisors of c,
then
\begin{displaymath}
\prod_{d|c}\;G_{d}\;=\;q_{1}^{\Theta(q_{1})} \ldots\;
q_{\omega}^{\Theta(q_{\omega})} \prod_{\substack{2\leq p<c,\\
(p,\;c)=1}}p^{\;\left[\frac{c}{p}\right]},
\end{displaymath}
where
\begin{displaymath}
\Theta(q_{i})=\sum_{\substack{d=0(q_{i}) \\
d|c}}\frac{\varphi(d)}{2}+\sum_{\substack{d \neq 0(q_{i}) \\
d|c}}E_{d}(q_{i}).
\end{displaymath}
\textbf{Proof.} Applying Theorem 1 to the divisors $d$ of $c$ one
has
\begin{displaymath}
G_{d}=R(d)^{\frac{\varphi d}{2}}\prod_{\substack{2\leq p<d\\
(p,\;d)=1}}p^{E_{d}(p)}.
\end{displaymath}
Multiplying over all $\tau (c)$ divisors and using the, easily
established, fact that
\begin{displaymath}
\prod_{\substack{x+y=c \\
x\leq y}}R(xyc)=R(c)^{\left[\frac{c}{2}\right]}\prod_{\substack{2\leq p<c\\
(p,\;c)=1}}p^{\;\left[\frac{c}{p}\right]},
\end{displaymath}
 gives for $\prod_{d|c}G_{d}$ the result.
\par
\vspace{20pt}
The corollary will not be used in the sequel.
\par
\vspace{20pt} We shall now prove certain Lemmas regarding the
function $E_{c}(x)$, and state, without proof, some well known
facts from the elementary theory of primes, so as not to interrupt
the main body of the proof. Absolute constants will be denoted by
$k_{i}$, indexed in the order they first appear.
\par
\textbf{Lemma 1.} For $x > 0$ (actually for any $x\neq 0)$
\begin{displaymath}
Max \hspace{2pt}\left( 0,\;\;\frac{\varphi(c)}{x}-2^{\omega
-1}\right)<E_{c}(x)<\frac{\varphi(c)}{x}+2^{\omega -1}.
\end{displaymath}
\par
\textbf{Proof.} By definition one has
\par
\begin{center}
$\displaystyle{\frac{1}{x}\frac{c}{1}-\binom{\omega}{0}<\left[\frac{1}{x}\frac{c}{1}\right]\leq
\frac{1}{x}\frac{c}{1}}$\\ \vspace{20pt}
$\displaystyle{-\sum_{1\leq
i\leq\omega}\frac{1}{x}\frac{c}{q_{i}} \leq -\sum_{1\leq
i\leq\omega}\left[\frac{1}{x}\frac{c}{q_{i}}\right]<\binom{\omega}{1}-\sum_{1\leq
i\leq\omega} \frac{1}{x}\frac{c}{q_{i}}} $\\ \vspace{20pt}
$\displaystyle{\sum_{1\leq i\; j
\leq\omega}\frac{1}{x}\frac{c}{q_{i}q_{j}}-\binom{\omega}{2} <
\sum_{1\leq
i,j\leq\omega}\left[\frac{1}{x}\frac{c}{q_{i}q_{j}}\right] \leq
\sum_{1\leq i,j\leq\omega} \frac{1}{x}\frac{c}{q_{i}q_{j}}} $\\
\vspace{20pt} $\cdots \hspace{1in} \cdots \hspace{1in} \cdots \hspace{0.1in}$\\
\vspace{20pt} $\displaystyle{\frac{1}{x}\frac{c}{q_{1}\cdots
q_{\omega}}-\binom{\omega}{\omega}<(-1)^{\omega}\left[\frac{1}{x}\frac{c}{q_{1}\cdots
q_{\omega}}\right] \leq \frac{1}{x}\frac{c}{q_{1}\cdots
q_{\omega}}, \hspace{0.5in} \omega \equiv 0(2)}$ \\
\vspace{20pt} $\displaystyle{-\frac{1}{x}\frac{c}{q_{1}\cdots
q_{\omega}}\leq -1)^{\omega}\left[\frac{1}{x}\frac{c}{q_{1}\cdots
q_{\omega}}\right]<\binom{\omega}{\omega}-\frac{1}{x}\frac{c}{q_{1}\cdots
q_{\omega}}, \hspace{0.5in} \omega \equiv 1(2)}.$
\end{center}
Adding term-wise above inequalities, we have
\begin{center}
$\displaystyle{\frac{1}{x}\left\{c-\sum_{1\leq i \leq \omega
}\frac{c}{q_{i}}+\sum_{1\leq i,j \leq \omega
}\frac{c}{q_{i}q_{j}}-\cdots\right\} -\sum_{\nu\equiv
0(2)}\binom{\omega}{\nu}<E_{c}(x)<}$ \\ \vspace{20pt}
$\displaystyle{\frac{1}{x}\left\{c-\sum_{1\leq i \leq
\omega}\frac{c}{q_{i}}+\sum_{1\leq i,j \leq
\omega}\frac{c}{q_{i}q_{j}}-\cdots\right\} +\sum_{\nu\equiv
1(2)}\binom{\omega}{\nu}}.$
\end{center}
Considering that
\begin{displaymath}
\sum_{\nu\equiv 0(2)}\binom{\omega}{\nu}=\sum_{\nu\equiv
1(2)}\binom{\omega}{\nu}=2^{\omega -1},
\end{displaymath}
and since $E_{c}(x)$ is the \# of numbers $n \leq \dfrac{c}{x}$,
$(n,c)=1$, i.e. always $\geq 0$, we have, as required,
\begin{displaymath}
Max \hspace{2pt}\left( 0,\;\;\frac{\varphi(c)}{x}-2^{\omega
-1}\right)<E_{c}(x)<\frac{\varphi(c)}{x}+2^{\omega -1}.
\end{displaymath}
\vspace{20pt}\textbf{Lemma 2.}
\begin{displaymath}
E_{c}(x)\hspace{10pt}>\hspace{10pt}\left\{\begin{matrix}\
\hspace{5pt} \dfrac{\varphi(c)}{x}-2^{\omega -1} &\hspace{30pt}
\text{for}& \hspace{0.5in}0<x<\dfrac{\varphi(c)}{2^{\omega
-1}} \\
 \hspace{0.2in}1  & \hspace{30pt}\text{for}&
\hspace{0.5in}\dfrac{\varphi(c)}{2^{\omega -1}} \leq x < c .
\end{matrix}\right.
\end{displaymath}
\par
\textbf{Proof.} The expression $\dfrac{\varphi(c)}{x}-2^{\omega
-1}$ is positive for all $0<x<\dfrac{\varphi(c)}{2^{\omega-1}}$.
For $x \geq \dfrac{\varphi(c)}{2^{\omega-1}}$ the lowest limit of
$E_{c}(x)$ is $1$, since $x<c$.
\par
\vspace{20pt}\textbf{Lemma 3.} If one of the prime factors
$q_{i}$, has exponent $\boxed{\alpha_{i}\geq 2}$ then
\begin{displaymath}
E_{c}(q_{i})=\frac{\varphi(c)}{q_{i}}.
\end{displaymath}
\par
\textbf{Proof.} Writing $q_{i}$ instead of $x$ in $E_{c}(x)$, and
renaming the running indexes, one has
\begin{displaymath}
E_{c}(q_{i})=\left[\frac{1}{q_{i}}\frac{c}{1}\right]-\sum_{1\leq j
\leq
\omega}\left[\frac{1}{q_{i}}\frac{c}{q_{j}}\right]+\sum_{1\leq
j,k \leq
\omega}\left[\frac{1}{q_{i}}\frac{c}{q_{j}q_{k}}\right]-\cdots
+(-1)^{\omega}\left[\frac{1}{q_{i}}\frac{c}{q_{1}\cdots
q_{\omega}}\right].
\end{displaymath}
Since by supposition the exponent of $q_{i}$ is $\geq 2,$ the
numbers within the integral part brackets are integers so that we
can skip the brackets. This gives
\begin{displaymath}
\frac{1}{q_{i}}\frac{c}{1}-\sum_{1\leq j \leq
\omega}\frac{1}{q_{i}}\frac{c}{q_{j}}+\sum_{1\leq j,k \leq
\omega}\frac{1}{q_{i}}\frac{c}{q_{j}q_{k}}-\cdots
+(-1)^{\omega}\frac{1}{q_{i}}\frac{c}{q_{1}\cdots
q_{\omega}}=\frac{\varphi(c)}{q_{i}},
\end{displaymath}
as stated.
\par
\vspace{20pt} \textbf{Lemma 4.} If one of the prime factors
$q_{i}$ has exponent $\boxed{\alpha_{i} = 1}$, then
\begin{displaymath}
E_{c}(q_{i})=\frac{\varphi(c)}{q_{i}-1}-E_{c/q_{i}}(q_{i}).
\end{displaymath}
\par
\textbf{Proof.} For convenience, putting $\overline{c}=c/q_{i}$,
renaming indexes as above, and writing on the left hand side the
terms of the $E_{c}(q_{i})$ function vertically, we have
following equalities
\begin{center}
$\displaystyle{\left[\frac{1}{q_{i}}\frac{c}{1}\right]=\left[\frac{1}{q_{i}}\frac{q_{i}\overline{c}}{1}
\right]} $\\ \vspace{20pt} $\displaystyle{-\sum_{1\leq
j\leq\omega}\left[\frac{1}{q_{i}}\frac{c}{q_{j}}\right]=-\left[\frac{1}{q_{i}}
\frac{q_{i}\overline{c}}{q_{i}}\right]-\sum_{2\leq
j\leq\omega}\left[\frac{1}{q_{i}}\frac{q_{i}\overline{c}}{q_{j}}\right]}$\\
\vspace{20pt} $\displaystyle{+\sum_{1\leq
j,\;k\leq\omega}\left[\frac{1}{q_{i}}\frac{c}{q_{j}q_{k}}\right]=+\sum_{2\leq
j\leq\omega}\left[\frac{1}{q_{i}}
\frac{q_{i}\overline{c}}{q_{i}q_{j}}\right]+\sum_{2\leq j,\; k
\leq\omega}\left[\frac{1}{q_{i}}\frac{q_{i}\overline{c}}{q_{j}q_{k}}\right]}$
\\ \vspace{20pt}$\displaystyle{\ldots \hspace{1in} \ldots\hspace{1in} \ldots}$
\\ \vspace{20pt}$\displaystyle{(-1)^{\omega -1} \sum_{1 \leq
j_{\,\nu}\leq
\omega}\left[\frac{1}{q_{i}}\frac{c}{q_{i_{1}}q_{j_{1}} \cdots
q_{j_{\omega -1}} }\right]= (-1)^{\omega -1} \sum_{1 \leq
j_{\,\nu}\leq
\omega}\left[\frac{1}{q_{i}}\frac{q_{i}\overline{c}}{q_{i}q_{i_{1}}q_{j_{1}}
\cdots q_{j_{\omega -1}} }\right] + }$\\  \vspace{20pt}
$\displaystyle{ \hspace{3.3in}\hspace{10pt} (-1)^{\omega
-1}\left[\frac{1}{q_{i}}\frac{q_{i}\overline{c}}{q_{2} \cdots
q_{\omega } }\right]  }$ \\  \vspace{20pt}
$\displaystyle{(-1)^{\omega}\left[\frac{1}{q_{i}}\frac{c}{q_{1}q_{2}
\cdots
q_{\omega}}\right]=(-1)^{\omega}\left[\frac{1}{q_{i}}\frac{q_{i}\overline{c}}{q_{i}q_{2}
\cdots q_{\omega}}\right]}.$
\end{center}
Adding above equalities, the sum of the left hand side terms is,
as stated, $E_{c}(q_{i})$. The right hand side is equal to
\begin{center}
$\displaystyle{\left\{-\left[\frac{1}{q_{i}}\frac{q_{i}\overline{c}}{q_{i}}\right]
+\sum_{2\leq
j\leq\omega}\left[\frac{1}{q_{i}}\frac{q_{i}\overline{c}}{q_{i}q_{j}}\right]
+ \ldots \right .}$ \\ \vspace{20pt}$\displaystyle{\left .+
(-1)^{\omega -1} \sum_{1 \leq j_{\,\nu}\leq \omega
}\left[\frac{1}{q_{i}}\frac{q_{i}\overline{c}}{q_{i}q_{i_{1}}q_{j_{1}}
\cdots q_{j_{\omega -2}}
}\right]+(-1)^{\omega}\left[\frac{1}{q_{i}}\frac{q_{i}\overline{c}}{q_{i}q_{2}
\cdots q_{\omega} }\right] \right \}+}$ \\ \vspace{20pt}
$\displaystyle{\left\{\left[\frac{1}{q_{i}}\frac{q_{i}\overline{c}}{1}\right]
- \sum_{2\leq
j\leq\omega}\left[\frac{1}{q_{i}}\frac{q_{i}\overline{c}}{q_{j}}\right]
+ \sum_{2\leq
j,\;k\leq\omega}\left[\frac{1}{q_{i}}\frac{q_{i}\overline{c}}{q_{j}q_{k}}\right]
 - \cdots
+(-1)^{\omega
-1}\left[\frac{1}{q_{i}}\frac{q_{i}\overline{c}}{q_{2} \cdots
q_{\omega} }\right] \right\}.}$
\end{center}
Cancelling $q_{i}$ within the brackets, the first brace is clearly
equal to $-E_{\overline{c}}(q_{i})$. In the second brace, as all
$q_{j},\; q_{k},\;\ldots\;, 2 \leq j,\; k \leq \omega,\; \ldots$
divide $\overline{c}$, the numbers within the integral part
brackets are integers. We therefore can write it
\begin{center}
$\displaystyle{\overline{c} \left( 1- \sum_{2\leq
j\leq\omega}\frac{1}{q_{j}}- \sum_{2\leq
j\;k\leq\omega}\frac{1}{q_{j}q_{k}} + \cdots +(-1)^{\omega
-1}\frac{1}{q_{2}\cdots q_{\omega}}\right)=
\varphi(\overline{c})  }$
\end{center}
Substituting the braces by their found values and taking into
account that
 $\varphi(\overline{c})=\varphi(c)/(q_{i}-1)$, since
$(c, \overline{c}) = 1$, we get
\begin{displaymath}
E_{c}(q_{i})=\frac{\varphi(c)}{q_{i}-1}-E_{c/q_{i}}(q_{i}),
\end{displaymath}
which proves the Lemma.
\par
\vspace{20pt} \textbf{Lemma 5.}
\begin{displaymath}
e^{k_{1}c} < \prod_{2\leq p \leq c}p \hspace{5pt}< e^{k_{2}c}.
\end{displaymath}
\par
\textbf{Proof.} This is the multiplicative form of Tchebycheff 's
estimate for  $\sum\limits_{2\leq p \leq c}\log p$, namely,
\begin{displaymath}
k_{1}c < \sum_{2\leq p \leq c}\log p < k_{2}c,
\end{displaymath}
for $c\geq 2$, with $k_{1}$ $k_{2}$ positive absolute constants.
\par
\vspace{20pt} \textbf{Lemma 6.}
\begin{displaymath}
e^{-k_{3}}c < \prod_{2\leq p \leq c} p^{\frac{1}{p}} < e^{k_{3}}c.
\end{displaymath}
\par
\textbf{Proof.} This is the multiplicative form of Merten's
estimate for $\sum\limits_{2\leq p \leq c}\frac{1}{p}\log p$,
namely,
\begin{displaymath}
\log c -k_{3} < \sum_{2\leq p \leq c} {\frac{1}{p}}\log p < \log
c+k_{3},
\end{displaymath}
for $c \geq 2$, with $k_{3}$ a positive absolute constant.
\par\vspace{20pt}
We now state a result which gives a lower bound for the geometric
mean $G_{c}^{\frac{2}{\varphi(c)}}$ in terms of the prime factors
$q_{i}(c)$ and their exponents $\alpha_{i}$. \vspace{10pt}
\par
\vspace{20pt} \textbf{Theorem 2.}
\begin{displaymath}
G_{c}^{\frac{2}{\varphi(c)}} > k_{4}\prod_{1\leq i\leq
\omega}q_{i}^{2\alpha_{i}-1-\frac{2}{\varphi(c)}E_{c}(q_{i})}\left(
\frac{q_{i}-1}{2}\right)^{2},
\end{displaymath}
where $k_{4}$ is a positive absolute constant.
\par
\textbf{Proof.} We transform the expression given for $G_{c}$ in
Theorem 1, as follows:
\begin{center}
$\displaystyle{G_{c}=\prod_{1 \leq i \leq
\omega}q_{i}^{\frac{\varphi(c)}{2}} \prod_{\substack{2\leq p < c \\
(p,\;c)=1}} q_{i}^{E_{c}(p)}  }$ \\ \vspace{20pt}
$\displaystyle{=\prod_{1 \leq i \leq
\omega}q_{i}^{\frac{\varphi(c)}{2}}\prod_{1 \leq i \leq
\omega}q_{i}^{-E_{c}(q_{i})}\prod_{2 \leq p < c}p^{E_{c}(p)} .}$
\end{center}
Joining the first two products into one and splitting the third
product as indicated, we have
\begin{center}
$\displaystyle{G_{c}=\prod_{1 \leq i \leq
\omega}q_{i}^{\frac{\varphi(c)}{2}-E_{c}(q_{i})} \prod_{2 \leq p
< \frac{\varphi(c)}{2^{\omega -1}}} p^{E_{c}(p)} \prod_{
\frac{\varphi(c)}{2^{\omega -1}}<p<c} p^{E_{c}(p)} .}$
\end{center}
Applying Lemma 2 to the second and third product and splitting
the products in an obvious way, we get successively
\begin{center}
$\displaystyle{G_{c} > \prod_{1 \leq i \leq
\omega}q_{i}^{\frac{\varphi(c)}{2}-E_{c}(q_{i})} \prod_{2 \leq p
< \frac{\varphi(c)}{2^{\omega -1}}}
p^{\frac{\varphi(c)}{p}-2^{\omega -1}} \prod_{
\frac{\varphi(c)}{2^{\omega -1}}<p<c} p }$ \\ \vspace{20pt}
$\displaystyle{>\prod_{1 \leq i \leq
\omega}q_{i}^{\frac{\varphi(c)}{2}-E_{c}(q_{i})}
\left\{\prod_{2\leq p < \frac{\varphi(c)}{2^{\omega -1}}}
p^{\frac{1}{p}}\right\}^{\varphi(c)} \left\{\prod_{2\leq p <
\frac{\varphi(c)}{2^{\omega -1}}} p\right\}^{-2^{\omega -1}}
\left\{\prod_{2\leq p < \frac{\varphi(c)}{2^{\omega -1}}}
p\right\}^{-1} \prod_{2 \leq p < c}p   }.$
\end{center}
Joining the third and the fourth product into one, we have
\begin{displaymath}
G_{c} > \prod_{1 \leq i \leq
\omega}q_{i}^{\frac{\varphi(c)}{2}-E_{c}(q_{i})}
\left\{\prod_{2\leq p < \frac{\varphi(c)}{2^{\omega -1}}}
p^{\frac{1}{p}}\right\}^{\varphi(c)} \left\{\prod_{2\leq p <
\frac{\varphi(c)}{2^{\omega -1}}} p\right\}^{-(2^{\omega -1}+1)}
\prod_{2 \leq p < c}p .
\end{displaymath}
Applying Lemma 6 to the second product, Lemma 5 to the third and
fourth product, we have
\begin{displaymath}
G_{c} >  \prod_{1 \leq i \leq
\omega}q_{i}^{\frac{\varphi(c)}{2}-E_{c}(q_{i})}
\left\{e^{-k_{3}}\frac{\varphi(c)}{2^{\omega
-1}}\right\}^{\varphi(c)}\left\{e^{k_{2}\frac{\varphi(c)}{2^{\omega
-1}}}\right\}^{-(2^{\omega -1}+1)}e^{k_{1}c} .
\end{displaymath}
Summing the exponents of $e$, we have
\begin{displaymath}
G_{c} > \prod_{1 \leq i \leq
\omega}q_{i}^{\frac{\varphi(c)}{2}-E_{c}(q_{i})}\left(
\frac{\varphi(c)}{2^{\omega -1} }\right)^{\varphi(c)}
e^{-k_{3}\varphi(c)+k_{1}c-k_{2}\varphi(c)-k_{2}\frac{\varphi(c)}{2^{\omega
-1}}}.
\end{displaymath}
Raising this inequality to the power $\frac{2}{\varphi(c)}$ we
get for the geometric mean $G_{c}^{\frac{2}{\varphi(c)}}$,
\begin{displaymath}\hspace{0.3in}
G_{c}^{\frac{2}{\varphi(c)}}> \prod_{1 \leq i \leq
\omega}q_{i}^{1-\frac{2}{\varphi(c)}E_{c}(q_{i})}
\left(\frac{\varphi(c)}{2^{\omega -1}}\right)^{2}
e^{-2k_{3}+2k_{1}\frac{c}{\varphi(c)}-2k_{2}\left(1+\frac{1}{2^{\omega
-1 }}\right)}.\hspace{0.8in} (3)
\end{displaymath}
Evaluating the second parenthesis, we have
\begin{displaymath}\hspace{0.35in}
\left ( \frac{\varphi(c)}{2^{\omega -1}} \right )^{2}=4\prod_{1
\leq i \leq \omega}q_{i}^{2\alpha_{i}-2}\left ( \frac{q_{i}-1}{2}
\right )^{2} .
\end{displaymath}

\par
On the other hand, since $\dfrac{c}{\varphi(c)} > 1$ and
$\left(1+\dfrac{1}{2^{\omega -1}}\right)\leq 2$, the entire
exponent of $e$ appearing in (3) is $> -2k_{3}+2k_{1}-4k_{2}$.
Setting $k_{4}=4e^{-2k_{3}+2k_{1}-4k_{2}}$, as a new absolute
constant and substituting in (3), we finally get
\begin{displaymath}
G_{c}^{\frac{2}{\varphi(c)}} > k_{4}\prod_{1\leq i\leq
\omega}q_{i}^{2\alpha_{i}-1-\frac{2}{\varphi(c)}E_{c}(q_{i})}\left(
\frac{q_{i}-1}{2}\right)^{2},
\end{displaymath}
which was to be proved.
\par
Following main theorem gives a lower bound for the geometric mean
$G_{c}^{\frac{2}{\varphi(c)}}$ in terms of $c$ and its radical.
\par
\vspace{20pt} \textbf{Theorem 3.} For any given $\varepsilon > 0$
\begin{displaymath}
G_{c}^{\frac{2}{\varphi(c)}}\;>\;k_{\varepsilon}\;R(c)^{1-\varepsilon}\;c^{2},
\end{displaymath}
where $k_{\varepsilon}$ is a positive absolute constant,
depending on $\varepsilon$.
\par
\textbf{Proof.}Denote the expressions within the product of
theorem 2 by $F(q_{i},\alpha_{i})$,\hspace{5pt}$1\leq q_{i} \leq
\omega$. The proof is in three steps. The first step gives a
lower bound of $F(q_{i},\alpha_{i})$ for $\alpha_{i}\geq 2$. The
second step gives a lower bound of $F(q_{i},\alpha_{i})$ for
$\alpha_{i}=1$. The third step combines these results to prove the
theorem.
\par
\underline{Step 1}. For $\boxed{\alpha_{i}\geq 2}$ we have by
Lemma 3
\begin{center}
\begin{tabular}{ccl}
$ F(q_{i},\alpha_{i})$ & $ = $ &
$q_{i}^{2\alpha_{i}-1-\frac{2}{\varphi
(c)}E_{c}(q_{i})}\left (\frac{q_{i}-1}{2}\right )^{2}$ \\
\hspace{10pt}& $ = $ & $q_{i}^{2\alpha_{i}-1-\frac{2}{q_{i}}}\left
(\frac{q_{i}-1}{2}\right )^{2}$ \\
\hspace{10pt} & $ = $ & $ q_{i}^{2\alpha_{i}+1}\left (
q_{i}^{\frac{1}{q_{i}}} \right )^{-2}\left (
\frac{1}{2}-\frac{1}{2q_{i}}\right )^{2}.$
\end{tabular}
\end{center}
Since for $q_{i}\geq 2$ the product $\left (
q_{i}^{\frac{1}{q_{i}}} \right )^{-2}\left (
\frac{1}{2}-\frac{1}{2q_{i}}\right )^{2}$, increasing
monotonically, tends to $1/4$ for $q_{i}\rightarrow \infty$, we
can write, for any given $\varepsilon > 0$
\begin{displaymath}
F(q_{i},\alpha_{i}) >
q_{i}^{2\alpha_{i}+1}q_{i}^{-\varepsilon}=q_{i}^{2\alpha_{i}+1-\varepsilon},
\end{displaymath}
for all \underline{primes greater than $N_{\varepsilon}$}, where
$N_{\varepsilon}$ is a number depending only on $\varepsilon$.
\par
For \underline{primes smaller than $N_{\varepsilon}$}, we can
write
\begin{displaymath}
F(q_{i},\alpha_{i}) \;> \;q_{i}^{2\alpha_{i}+1}\left (
q_{i}^{\frac{1}{q_{i}}} \right )^{-2}\left (
\frac{1}{2}-\frac{1}{2q_{i}}\right )^{2}q_{i}^{-\varepsilon} \;
>\; \frac{1}{32}\;\;q_{i}^{2\alpha_{i}+1-\varepsilon},
\end{displaymath}
since, as said before, the product $\left (
q_{i}^{\frac{1}{q_{i}}} \right )^{-2}\left (
\frac{1}{2}-\frac{1}{2q_{i}}\right )^{2}$ is monotonically
increasing, and therefore its minimum is at $q_{i}=2$, i.e. is
equal to \\  $\left( 2^{\frac{1}{2}} \right )^{-2} \left(
\frac{1}{2}-\frac{1}{4}\right )^{2}=\frac{1}{32}$.
\par
$N_{\varepsilon}$ can be calculated. It is the abscissa where the
two curves $y =\left ( x^{\frac{1}{x}} \right )^{-2}\left (
\frac{1}{2}-\frac{1}{2x}\right )^{2}$ and $y=x^{-\varepsilon}$
cut each other. Consequently it is the positive root of the
equation
\begin{displaymath}
x-1=2x^{\frac{1}{x}+\frac{2-\varepsilon}{2}},
\end{displaymath}
which results after setting
\begin{displaymath}
\left ( x^{\frac{1}{x}} \right )^{-2}\left (
\frac{1}{2}-\frac{1}{2x}\right )^{2}=x^{-\varepsilon}.
\end{displaymath}
By elementary analysis $N_{\varepsilon}\rightarrow \infty$ for
$\varepsilon\rightarrow 0$.
\par
\vspace{20pt} \underline{ Step 2}. For $\boxed{\alpha_{i}=1}$ we
have by Lemma 4 and Lemma 1
\begin{center}
\begin{tabular}{ccl}
$ F(q_{i},\alpha_{i})$ & $ = $ &
$q_{i}^{2\alpha_{i}-1-\frac{2}{\varphi
(c)}E_{c}(q_{i})}\left (\frac{q_{i}-1}{2}\right )^{2}$ \\
\hspace{10pt}& $ = $ &
$q_{i}^{2\alpha_{i}-1-\frac{2}{\varphi(c)}\left [ \frac{\varphi
(c)}{q_{i}-1} - E_{c/q_{i}}(q_{i}) \right ]}\left (
\frac{q_{i}-1}{2} \right )^{2}$ \\
\hspace{10pt} & $ = $ &
$q_{i}^{2\alpha_{i}-1-\frac{2}{q_{i}-1}+\frac{2}{\varphi
(c)}E_{c/q_{i}}(q_{i})}\left ( \frac{q_{i}-1}{2} \right )^{2} $ \\
\hspace{10pt} & $ > $ &
$q_{i}^{2\alpha_{i}-1-\frac{2}{q_{i}-1}}\left ( \frac{q_{i}-1}{2}
\right )^{2}=\; q_{i}^{2\alpha_{i}+1}\left
(q_{i}^{\frac{1}{q_{i}-1}}\right )^{-2} \left
(\frac{1}{2}-\frac{1}{2q_{i}}\right )^{2}. $
\end{tabular}
\end{center}
Since for $q_{i}\geq 2$ the product $\left
(q_{i}^{\frac{1}{q_{i}-1}}\right )^{-2} \left
(\frac{1}{2}-\frac{1}{2q_{i}}\right )^{2}$, increasing
monotonically, tends to $\frac{1}{4}$ for $q_{i}\rightarrow
\infty$, we can write
\begin{displaymath}
F(q_{i},\alpha_{i}) >\;
q_{i}^{2\alpha_{i}+1}q_{i}^{-\varepsilon}=q_{i}^{3-\varepsilon},
\end{displaymath}
for all \underline{primes greater than $M_{\varepsilon}$}, where
$M_{\varepsilon}$ is a certain number depending only on
$\varepsilon$.
\par
For \underline{primes smaller than $M_{\varepsilon}$} we can write
\begin{center}
\begin{tabular}{ccl}
$ F(q_{i},\alpha_{i})$ & $ > $ & $q_{i}^{2\alpha_{i}+1}\left
(q_{i}^{\frac{1}{q_{i}-1}}\right )^{-2} \left
(\frac{1}{2}-\frac{1}{2q_{i}}\right )^{2}q_{i}^{-\varepsilon}$ \\
\hspace{10pt} & $ > $ &
$\frac{1}{64}q_{i}^{2\alpha_{i}+1-\varepsilon}=\frac{1}{64}q_{i}^{3-\varepsilon},$
\end{tabular}
\end{center}
since the product $\left ( q_{i}^{\frac{1}{q_{i}-1}}\right
)^{-2}\left ( \frac{1}{2}-\frac{1}{2q_{i}}\right )^{2}$ is
monotonically increasing, and therefore its minimum is at
$q_{i}=2$, i.e. is equal to $\left ( 2^{-2} \right ) \left (
\frac{1}{2}-\frac{1}{4} \right )^{2}=\frac{1}{64}$.
\par
$M_{\varepsilon}$ can be calculated as a function of
$\varepsilon$. It is the abscissa where the two curves $ y=\left
(x^{\frac{1}{x-1}}\right )^{-2} \left
(\frac{1}{2}-\frac{1}{2x}\right )^{2}$ and $y=x^{-\varepsilon}$
cut each other. It is therefore the positive root of the equation
\begin{displaymath}
x-1=2x^{\frac{1}{x-1}+\frac{2-\varepsilon}{2}},
\end{displaymath}
which results after setting
\begin{displaymath}
\left (x^{\frac{1}{x-1}}\right )^{-2}\left
(\frac{1}{2}-\frac{1}{2x}\right )^{2} = x^{-\varepsilon}.
\end{displaymath}
By elementary analysis $M_{\varepsilon}\rightarrow \infty$ for
$\varepsilon\rightarrow 0$.
\par
\vspace{20pt} \underline{Step 3}. Let $p$ run through the primes
$q_{i}$, $1\leq i \leq \omega$ of $c$, and denote by $D(p)$ the
corresponding exponents $\alpha_{i}$.
\par
All the integers $q_{i}^{\alpha_{i}}$, $1\leq i \leq \omega$, fall
into one of the following, mutually exclusive, classes specified
under Steps 1 and Step 2, namely
\par
\begin{center}
\begin{tabular}{ccl}
$2\leq p\leq \left [ N_{\varepsilon} \right ]$ &\; $D(p) \geq 2$
\;& $\omega_{1}(\varepsilon)$ \\ \\
$\left [ N_{\varepsilon} \right ]+1 \leq p $&\; $D(p) \geq 2 $
\;&$ \omega_{2}(\varepsilon,c)$ \\ \\
$2\leq p \leq \left [ M_{\varepsilon} \right ] $& \;$ D(p)=1 $
\;& $\omega_{3}(\varepsilon) $\\ \\
$\left [ M_{\varepsilon} \right ]+1 \leq p $& \;$ D(p)=1$ \;&
$\omega_{4}(\varepsilon,c).$
\end{tabular}
\end{center}
\par
It is clear from their definition that
\begin{displaymath}
\omega_{1}(\varepsilon)\;<\;\pi(N_{\varepsilon})
\end{displaymath}
\begin{displaymath}
\omega_{3}(\varepsilon)\;<\;\pi(M_{\varepsilon})
\end{displaymath}
where $\pi(x)$ denotes the number of primes not exceeding $x$.
\par
As indicated $\omega_{2}$ and $\omega_{4}$ depend on both
$\varepsilon$ and $c$.
\par
 \vspace{20pt} Accordingly, the product $ \prod\limits_{1\leq
i\leq\omega}F(q_{i},\alpha_{i})$ extended over above classes can
be written
\begin{displaymath}
\prod_{1\leq i \leq\omega}F(q_{i},\alpha_{i})= \prod_{\substack{2\leq p < [N_{\varepsilon}] \\
D(p)\geq 2}}F(p,D(p))
\prod_{\substack{[N_{\varepsilon}]+1\leq p  \\
D(p)\geq 2}}F(p,D(p)) \\
\end{displaymath}
\begin{displaymath}
\hspace{5cm}
\prod_{\substack{2\leq p \leq [M_{\varepsilon}] \\
D(p)=1}}F(p,D(p))
\prod_{\substack{[M_{\varepsilon}]+1\leq p \\
D(p)=1}}F(p,D(p)).
\end{displaymath}
Applying to these four products the lower bounds found in Step 1
and Step 2 for the respective expressions $F(q_{i},\alpha_{i})$,
we obtain
\begin{displaymath}
\prod_{1\leq i \leq\omega}F(q_{i},\alpha_{i})>
\frac{1}{32^{\omega_{1}(\varepsilon)}}
\prod_{\substack{2\leq p < [N_{\varepsilon}] \\
D(p)\geq 2}}p^{2D(p)+1-\varepsilon}
\prod_{\substack{[N_{\varepsilon}]+1\leq p  \\
D(p)\geq 2}}p^{2D(p)+1-\varepsilon} \\
\end{displaymath}
\begin{displaymath}
\hspace{4cm}\frac{1}{64^{\omega_{3}(\varepsilon)}}
\prod_{\substack{2\leq p \leq [M_{\varepsilon}] \\
D(p)=1}}p^{2D(p)+1-\varepsilon}
\prod_{\substack{[M_{\varepsilon}]+1\leq p \\
D(p)=1}}p^{2D(p)+1-\varepsilon}.
\end{displaymath}
Since, as said above, the four classes cover, by definition, the
whole range of the integers $q_{i}^{\alpha_{i}}$, $1\leq i \leq
\omega$, this inequality can be written
\begin{displaymath}
\prod_{1\leq i \leq\omega}F(q_{i},\alpha_{i})>
\frac{1}{2^{5\omega_{1}(\varepsilon)+6\omega_{3}(\varepsilon)}}
\prod_{1\leq i \leq \omega}q_{i}^{2\alpha_{i}+1-\varepsilon}=
\frac{1}{2^{5\omega_{1}(\varepsilon)+6\omega_{3}(\varepsilon)}}\;R(c)^{1-\varepsilon}\;
c^{2} .
\end{displaymath}
By Theorem 2 we therefore have
\begin{displaymath}
G_{c}^{\frac{2}{\varphi(c)}} = k_{4}\prod_{1\leq i
\leq\omega}F(q_{i},\alpha_{i})>\frac{k_{4}}{2^{5\omega_{1}(\varepsilon)+6\omega_{3}(\varepsilon)}}
\;R(c) ^{1-\varepsilon} \;c^{2}
\end{displaymath}
\begin{displaymath}
\hspace{4.8cm} >\;
\frac{k_{4}}{2^{5\pi(N_{\varepsilon})+6\pi(M_{\varepsilon})}}
\;R(c) ^{1-\varepsilon} \;c^{2}
\end{displaymath}
Setting
$\frac{k_{4}}{2^{5\pi(N_{\varepsilon})+6\pi(M_{\varepsilon})}}$
as an absolute constant depending on $\varepsilon$, we finally get
\begin{displaymath}
G_{c}^{\frac{2}{\varphi(c)}}>
k_{\varepsilon}\;R(c)^{1-\varepsilon}\;c^{2},
\end{displaymath}
as claimed by Theorem 3.
\par
\vspace{20pt}
From above theorem results immediately
\par
\textbf{Theorem 4.} For any positive integer there are partitions
$c=a+b$, with positive coprime integers $a$ and $b$, such that
\begin{displaymath}
k_{\varepsilon}\;R(c)^{1-\varepsilon}\;c^{2} < R(abc).
\end{displaymath}
\textbf{Proof.} Because of Theorem 3, all radicals of equations
(1) which are greater than $G_{c}^{\frac{2}{\varphi(c)}}$ satisfy
a fortiori the condition of Theorem 4.
\par
For said partitions this is a substantially stronger result than
the abc-conjecture
$$\kappa_{\varepsilon}\;c^{\frac{1}{1+\varepsilon}}\hspace{5pt}\;
<\; \hspace{5pt} R(abc), \hspace{15pt} \varepsilon > \;0.$$
\par
\vspace{20pt} We shall now obtain an upper bound for
$G_{c}^{\frac{2}{\varphi(c)}}$. This will depend on $\omega$, the
number of prime factors of $c$. To this end we prove following
\par
\textbf{Theorem 5.}
$$G_{c}^{\frac{2}{\varphi(c)}} <\;
k_{5}\;k_{6}^{3^{\omega}}\;R(c)\;c^{2},$$ where $k_{5}>0$ and
$k_{6}>1$ are absolute constants.
\par
\textbf{Proof.} Applying Lemma 1, Lemma 5 and Lemma 6, we get
successively
\begin{center}
\begin{tabular}{clll}
$G_{c}$ & $=$ & $R(c)^{\frac{\varphi(c)}{2}}\prod\limits_{\substack{2\leq p < c \\
(p,\;c)=1}}p^{E_{c}(p)}$ & by Theorem 1 \\
\hspace{10pt}&$ = $ & $\prod\limits_{1\leq i\leq
\omega}q_{i}^{\frac{\varphi(c)}{2}-E_{c}(q_{i})}\prod\limits_{2\leq
p\leq
c}p^{E_{c}(p)}$ & \\
\hspace{10pt} &$<$& $\prod\limits_{1\leq i\leq
\omega}q_{i}^{\frac{\varphi(c)}{2}-E_{c}(q_{i})}\prod\limits_{2\leq
p\leq
c}p^{\frac{\varphi(c)}{p}+2^{\omega -1}}$& by Lemma 1 \\
\hspace{10pt}& $<$ & $\prod\limits_{1\leq i\leq
\omega}q_{i}^{\frac{\varphi(c)}{2}-E_{c}(q_{i})}
\left\{\prod\limits_{2\leq p \leq c}p^{\frac{1}{p}}
\right\}^{\varphi(c)}
\left\{\prod\limits_{2\leq p \leq c}p\right\}^{2^{\omega -1}}$ & \\
\hspace{10pt} & $<$ & $ \prod\limits_{1\leq i\leq
\omega}q_{i}^{\frac{\varphi(c)}{2}-E_{c}(q_{i})}\left(e^{k_{3}}c\right)
^{\varphi(c)}\left(e^{k_{2}c}\right) ^{2^{\omega -1}}$ & by Lemma
5, 6 \\
\hspace{10pt} & $<$ & $e^{k_{3}\varphi(c)+k_{2}2^{\omega
-1}c}\;\;c^{\varphi(c)}\prod\limits_{1\leq i\leq
\omega}q_{i}^{\frac{\varphi(c)}{2}-E_{c}(q_{i})}$ & .
\end{tabular}
\end{center}

 Raising the
last inequality to the power $\frac{2}{\varphi(c)}$ we have
$$G_{c}^{\frac{2}{\varphi(c)}}\;<\;e^{2k_{3}+k_{2}2^{\omega}\frac{c}{\varphi(c)}}\;\;c^{2} \prod_{1\leq i\leq
\omega}q_{i}^{1-\frac{2}{\varphi(c)}E_{c}(q_{i})}.$$ Since
$\dfrac{2q_{i}}{q_{i}-1} \leq 3\;$ for all $\;q_{i}\geq 3\;$ and
 $\dfrac{2q_{i}}{q_{i}-1}=4\;$ for $\;q_{i}=2$, we have
$$ 2^{\omega}\frac{c}{\varphi(c)}=\prod_{1\leq i\leq
\omega}\frac{2q_{i}}{q_{i}-1}\;<\;2\cdot 3^{\omega}.$$
\vspace{20pt}
 Hence
$$G_{c}^{\frac{2}{\varphi(c)}}\;<\;e^{2k_{3}+2k_{2}3^{\omega}}\;c^{2}\;\prod_{1\leq i\leq
\omega}q_{i}^{1-\frac{2}{\varphi(c)}E_{c}(q_{i})}$$
$$<\;e^{2k_{3}+2k_{2}3^{\omega}}\;c^{2}\;\prod_{1\leq i\leq
\omega}q_{i}.\hspace{0.18in}$$
Setting $k_{5}=e^{2k_{3}}$ and
$k_{6}=e^{2k_{2}}$ as absolute constants, we obtain
$$G_{c}^{\frac{2}{\varphi(c)}}\;<\;k_{5}\;k_{6}^{3^{\omega}}\;R(c)\;c^{2},$$
as required.
\par
\vspace{20pt} \textbf{Remark.} Combining Theorem 4 and Theorem 5
we have
$$k_{\varepsilon}\;R(c)
^{1-\varepsilon}\;c^{2}\;<\;G_{c}^{\frac{2}{\varphi(c)}}\;<\;k_{5}\;k_{6}^{3^{\omega}}\;R(c)\;c^{2},$$
and dividing by $R(c)c^{2}$ we get
$$k_{\varepsilon}\;R(c)
^{-\varepsilon}\;<\;\frac{G_{c}^{\frac{2}{\varphi(c)}}}{R(c)c^{2}}\;<\;k_{5}\;k_{6}^{3^{\omega}}.$$
Letting $c$ run through the numbers $c_{x}=q_{1}^{x_{1}}\cdots
q_{\omega}^{x_{\omega}},\hspace{3pt}1\leq x_{i} <\infty$, which
all have the same radical $R(c)=R$ we have
\par
\begin{displaymath}
k_{\varepsilon}\;R^{-\varepsilon}\;\leq\; \liminf\limits_{1\leq
x_{i} < \infty}\;\;
\frac{G_{c_{x}}^{\frac{2}{\varphi(c_{x})}}}{R\;c_{x}^{2}} \;\leq\;
\limsup\limits_{1\leq x_{i} <\infty}\;\;
\frac{G_{c_{x}}^{\frac{2}{\varphi(c_{x})}}}{R\;c_{x}^{2}}\;
 \;\leq\; k_{5}\;k_{6}^{3^{\omega}}.
\end{displaymath}
\par

\vspace{20pt} \textbf{Acknowledgment.} I am indebted to Peter
Krikelis, Department of Mathematics, University of Athens, for his
unfailing assistance.

\end{document}